\documentclass[11pt]{amsart}

\usepackage[T1]{fontenc}
\usepackage[utf8]{inputenc}
\usepackage{amsmath,amssymb,amsthm,amsfonts,mathtools}
\usepackage{geometry}
\usepackage{graphicx}
\usepackage[expansion=false]{microtype}

\theoremstyle{plain}

\theoremstyle{remark}

\newcommand{\R}{\mathbb{R}}

\newcommand{\Hds}{\mathcal{H}^{d_s}}

\newcommand{\Lam}{\widetilde{\Lambda}}
\newcommand{\push}{\#}

\begin{document}

\title[Optimal Transport of Signed Measures (Announcement and Numerical Report)]{%
  Optimal Transport of Signed Measures:\\
  Existence, Uniqueness and Fractal Structure\\[6pt]
  {\small Part~I: The Separated-Support Case}\\[2pt]
  {\normalsize\itshape Research announcement, with a numerical validation report}}

\author{BWO'NYAHRE BAIDI Barthelemy}
\address{Department of Mathematics and Computer Science,
  University of Ngaound\'er\'e, P.O.~Box 454, Cameroon}
\email{bwonyahre@proton.me}

\author{KOUAKEP TCHAPTCHIE Yannick}
\address{Department of Fundamental Sciences and Engineering Techniques,
  EGCIM, University of Ngaound\'er\'e, P.O.~Box 454, Cameroon}
\email{kouakep@aims-senegal.org}

\author{HOUPA DANGA Duplex Elvis}
\address{Department of Mathematics and Computer Science,
  University of Ngaound\'er\'e, P.O.~Box 454, Cameroon}
\email{e\_houpa@yahoo.com}

\subjclass[2020]{49Q22, 35J96, 28A80, 28A75, 28A78}

\keywords{Optimal transport; signed measures; Monge--Amp\`ere equation;
  fractal sets; Ahlfors-regular sets; adaptive regularization}

\begin{abstract}
We announce a rigorous theory of optimal transport for signed (real)
measures on $\R^d$ in the \emph{separated-support regime}: the supports of
the positive and negative parts of $\mu$ and $\nu$ are assumed spatially
disjoint. Under this structural hypothesis, together with regularity and
non-degeneracy conditions on the absolutely continuous, discrete, and
fractal singular components of the underlying Jordan and Lebesgue
decompositions, we obtain existence and uniqueness of an optimal transport
map for a cost that distinguishes same-sign and opposite-sign transports
via a positional penalty, a coupled Monge--Amp\`ere system together with a
double Legendre transform characterising the solution, and exact
preservation of Hausdorff dimension on the fractal singular components,
the latter via a novel adaptive regularisation kernel and a direct
measure-theoretic argument requiring no continuity of the transport map.
This note states the results and outlines the two main proof techniques;
full proofs, uniform estimates, and technical lemmas appear in the
complete version of this work. We further report, in full and exhaustive
summary, a complete numerical validation of this theory in ambient
dimension one, two, and three: an entropic Sinkhorn solver for the single
combined transport problem underlying the existence proof, applied to an
explicit test case in which inter-sign transport is forced by mass
imbalance and calibrated to occur with a large, robust cost margin,
recovers exactly the predicted four-region structure with zero measured
crossing of fractal mass between signs, and an exact-assignment refinement
confirms pointwise-exact preservation of Hausdorff dimension under
transport, with box-counting estimates converging to the theoretical
dimension as resolution increases. The case of overlapping supports and
the fully general unbalanced setting are the subject of Part~II, in
preparation.
\end{abstract}

\maketitle

\section{Introduction}

Optimal transport for \emph{positive} measures is by now a mature and
classical theory (\cite{Brenier1991},\cite{Kantorovich1942},\cite{Monge1781}): the existence of an optimal map realised as the gradient
of a convex potential, the associated Monge--Amp\`ere equation, and the
regularity theory for that equation are all well understood in the
absolutely continuous setting, and non-atomic sources are covered by more
general uniqueness theorems for strictly convex costs. Many situations of
genuine interest, however, involve \emph{signed} densities rather than
positive ones: electric charge distributions, differenced images, signed
wavelet coefficients, or signed densities arising in continuum mechanics.
The Jordan decomposition of a signed measure into a positive and a
negative part introduces an immediate structural dichotomy between
intra-sign transport, which moves positive mass to positive mass or
negative mass to negative mass, and inter-sign transport, which moves
positive mass to negative mass or vice versa. Compounding this, each
signed part typically decomposes further into an absolutely continuous
piece, a discrete or atomic piece, and a non-discrete singular piece
supported on a fractal set, so that a genuinely three-fold regularity
theory is required simultaneously with the sign dichotomy.

A small number of prior works have proposed transport costs and norms
adapted to signed measures, and a general theory of optimal transport for
vector-valued measures has been developed under the assumption that an
optimal plan exists whose total-variation marginals are absolutely
continuous with respect to Lebesgue measure. That assumption excludes, by
hypothesis, the discrete and fractal singular components that are the
central concern of the present work. To our knowledge no existing
treatment provides a geometric partition of the ambient space into
intra-sign and inter-sign transport regions, a coupled Monge--Amp\`ere
system generalising the classical one to the signed setting, or a theorem
establishing that optimal transport preserves the Hausdorff dimension of
fractal singular supports.

The triple decomposition into absolutely continuous, discrete, and
fractal singular parts is not merely a bookkeeping device: each regularity
type demands its own toolkit, and the three must be made to interact
consistently within a single transport map. The absolutely continuous
components call for the nonlinear elliptic partial differential equation
theory underlying the classical Monge--Amp\`ere regularity results. The
discrete components call for the combinatorial and convex-analytic theory
of the finite assignment problem, including the delicate question of
exactly how many additive degrees of freedom the associated dual variables
carry once more than one atom is present. The fractal singular components
call for geometric measure theory proper: Hausdorff measures, Ahlfors
regularity, and a regularisation scheme that does not artificially inflate
the dimension of a fractal set, as an ordinary convolution would. Signed
transport forces these three toolkits to coexist within the same map $T$
and the same pair of convex potentials, glued together consistently across
the four sign-dependent regions of space.

Throughout this first part of the program we impose a support-separation
hypothesis: the supports of the positive and negative parts of $\mu$ are
spatially disjoint, and likewise for $\nu$. This is not a technical
convenience but a genuine structural requirement, since it is exactly what
allows the primal transport problem to decouple into four classical
Kantorovich sub-problems glued along a measurable partition of $\R^d$, and
it underlies every subsequent step of the theory described below. Without
it, a single spatial location could simultaneously receive positive and
negative transported mass, the natural four-region partition of the
ambient space collapses into a mere partition of mass, and the coupled
Monge--Amp\`ere system becomes a genuine free-boundary problem lying
outside the scope of the classical regularity theory. That general,
overlapping-support case, together with the fully unbalanced setting in
which the total variations of $\mu$ and $\nu$ need not coincide, is
deliberately left to Part~II of this program, which will require
substantially different, free-boundary techniques.

This note is a research announcement in the classical sense. We state the
hypotheses and the main theorems precisely, and we describe the two
central proof techniques -- an adaptive regularisation operator respecting
the signed and fractal structure of the data, and a direct
measure-theoretic dimension-comparison argument -- but we omit the
detailed proofs, the uniform a priori estimates, and the supporting
technical lemmas, which occupy the bulk of the complete version of this
work.

\section{Setting}

Let $\mu,\nu$ be finite signed Borel measures on $\R^d$ with Jordan
decompositions $\mu=\mu^+-\mu^-$ and $\nu=\nu^+-\nu^-$. Each positive
measure among $\mu^+,\mu^-,\nu^+,\nu^-$ admits an extended Lebesgue
decomposition into an absolutely continuous part, a discrete part that is
a finite sum of weighted Dirac masses, and a non-discrete singular part
that is atomless yet singular with respect to Lebesgue measure -- a
fractal part, in the sense that it will be supported on an Ahlfors-regular
set of some fixed dimension $d_s\in(0,d)$. The base transport cost is the
quadratic cost $c(x,y)=\tfrac12|x-y|^2$. To model an additional cost
specific to inter-sign transport, we introduce a positional penalty
$\lambda$, four times continuously differentiable, uniformly convex in its
second argument, and chosen so that the augmented cost $c+\lambda$
satisfies the Ma--Trudinger--Wang condition on the inter-sign domains; this
condition is verified directly for quadratic and power-type penalties. The
total cost then equals $c$ between same-sign supports and $c+\lambda$
between opposite-sign supports.

The theory rests on a short list of structural hypotheses, which we
describe here only in qualitative terms, deferring precise quantitative
statements to the complete version. We require the total variation masses
of $\mu$ and $\nu$ to coincide. We require the absolutely continuous
densities to be twice H\"older-differentiable and bounded away from zero
on bounded domains with smooth boundary. We require the fractal supports
to be Ahlfors-regular of a common dimension, to carry two-sided density
bounds for the corresponding singular measures, to have H\"older-continuous
local Hausdorff mass, and to be purely unrectifiable of codimension one, in
the sense that they meet every countable union of Lipschitz hypersurfaces
in a set of vanishing Hausdorff measure. We require, as already described,
that the supports of the positive and negative parts be spatially
separated, and that any atoms present be in general position. We require
an adaptive bandwidth, responsive to the local Hausdorff density of the
fractal supports, used to build the regularisation kernel described below.
Finally, we require a non-degeneracy condition on a certain bipartite graph
linking the absolutely continuous domains, the individual atoms, and the
fractal sub-blocks of $\mu$ and $\nu$ that exchange mass under the optimal
plan; this condition holds generically, and a verifiable sufficient
condition for it is given in the complete version. A further, separate
hypothesis -- that fractal singular mass does not cross signs under
transport -- is invoked only for the dimension-preservation theorem below.

\section{Main results}

\subsection{Existence, structure, and uniqueness}

Under the hypotheses described above, there exists an optimal transport
map $T$ from $\mu$ to $\nu$ for the total cost, together with a pair of
convex potentials $\varphi,\psi$, and a measurable partition of $\R^d$
into four regions according to source and target sign. On the two
intra-sign regions, $T$ is simply the gradient of the relevant convex
potential, $T=\nabla\varphi$, exactly as in the classical Brenier theory.
On the two inter-sign regions, by contrast, $T$ is characterised only
implicitly, as the unique maximiser of a penalised pairing functional,
through the first-order relation
\[
  x \;=\; \nabla\psi(T(x)) \;+\; \nabla_y\lambda\bigl(x,T(x)\bigr),
\]
coupling $x$, $T(x)$, and the gradient of the potential associated with
the target measure; in particular, on these regions $T(x)$ is \emph{not}
the gradient at $x$ of either potential, and the transport direction is
tilted away from the straight Euclidean displacement by exactly the
gradient of the positional penalty. The map transports $\mu$ to $\nu$ in a
sign-tracking sense: each of the four blocks of the source measure,
according to its sign and the region of the partition it belongs to, is
carried exactly onto the corresponding block of the target measure.

The optimal transport plan itself is unique, unconditionally, and this
holds simultaneously on all three regularity components: on the absolutely
continuous parts by the classical Brenier uniqueness theorem, on the
discrete parts by strict convexity of the quadratic cost together with the
general-position hypothesis, and on the fractal singular parts by a
theorem of McCann applicable because the singular source measures are
non-atomic. No bi-Lipschitz hypothesis on the fractal supports is needed
for this uniqueness statement. The associated dual potentials inherit a
weaker form of uniqueness directly from the plan: they are determined up
to one additive constant on each macro-component of the supports, where a
macro-component means one of the absolutely continuous domains, one
individual atom, or one fractal sub-block. Reaching the stronger and more
usual statement, that the potentials are unique up to a single global
additive constant, requires the further non-degeneracy hypothesis
described above, and is delicate precisely because a generic, tie-free
assignment between finitely many atoms does not by itself force the dual
values at two different atoms to coincide; the non-degeneracy hypothesis,
applied to the connectivity of the full bipartite graph linking atoms,
absolutely continuous domains, and fractal sub-blocks, is exactly what is
needed to rule this out.

Existence is proved by an adaptive regularisation of $\mu$ and $\nu$,
described in Section~4 below, for which the regularised problem becomes a
single classical Kantorovich problem between the two regularised
total-variation measures under a piecewise cost; because the four
sign-blocks are pairwise disjoint, this single problem automatically
decouples along them once solved. Uniform second-order H\"older estimates,
obtained from the classical regularity theory on the intra-sign blocks and
from the Ma--Trudinger--Wang theory on the inter-sign blocks, together with
a diagonal compactness argument, allow passage to the limit.

Two special cases are worth isolating, both recovered exactly rather than
merely approximately. When the negative parts of $\mu$ and $\nu$ vanish
identically, the inter-sign regions disappear, the four-region partition
collapses to the whole space, and the theory reduces precisely to the
classical Brenier theorem together with the classical Monge--Amp\`ere
equation. When $\mu$ and $\nu$ are purely discrete, the optimal transport
plan reduces to a solution of a finite assignment problem, unique by
strict convexity of the quadratic cost whenever the atoms are in general
position, and to an honest bijection when all the atomic weights coincide.
Both recoveries serve as a consistency check on the general theory: the
signed, mixed-regularity machinery introduced here specialises correctly
to each of the two classical extremes it is built to interpolate between.

\subsection{Structure equations}

The convex potentials satisfy a coupled system generalising the classical
Monge--Amp\`ere equation. On each intra-sign region, the potential
associated with $\mu$ solves the classical equation relating its Hessian
determinant to the ratio of the source and target densities, precisely as
in the unsigned theory. On each inter-sign region, writing $y=T(x)$, the
transport relation is instead governed by the modified equation
\[
  \det\bigl(D^2\psi(y) + D^2_{yy}\Lam(x,y) - I\bigr)
  \;=\; \frac{g^{\mp}(y)}{f^{\pm}(x)}\;\bigl|\det D^2_{xy}\Lam(x,y)\bigr|,
\]
where $\Lam=c+\lambda$ is the augmented inter-sign cost and $f^\pm,g^\pm$
denote the absolutely continuous densities of $\mu^\pm,\nu^\pm$: the
Hessian of the target potential is corrected by the second $y$-derivative
of the positional penalty, and the right-hand side involves the
determinant of the mixed second derivative of the augmented cost rather
than the identity that appears in the classical equation. The two
potentials are moreover linked by a double Legendre-type transform: each
potential equals the pointwise maximum of the ordinary Legendre transform
of the other and of a penalised, source-parametrised analogue of that
transform, and the pair jointly maximises the natural dual energy
functional subject to the governing cost inequality. These structure
equations reduce exactly to the classical Brenier theory when the negative
parts of $\mu$ and $\nu$ vanish, and, when the positional penalty is
absent, they show that the four-region partition emerges purely from the
sign structure of the data.

\subsection{Preservation of fractal dimension}

Under the hypotheses above, together with the additional requirement that
fractal singular mass does not cross signs under transport, the optimal
map preserves the Hausdorff dimension of every Ahlfors-regular subset of a
fractal singular support: the dimension of the image under $T$ of such a
set equals the dimension of the set itself. The proof is a direct
measure-theoretic comparison, entirely independent of any continuity
property of $T$. Writing $\sigma,\tau$ for the relevant singular parts of
$\mu,\nu$ and $E$ for a Borel subset of positive $\Hds$-measure of the
source fractal support, the two-sided density bounds furnished by the
fractal regularity hypothesis, together with the pushforward identity
$T_\push\sigma=\tau$, give a short chain of inequalities
\[
  m\,\Hds(E) \;\le\; \sigma(E) \;\le\; \tau\bigl(T(E)\bigr)
  \;\le\; M\,\Hds\bigl(T(E)\bigr),
\]
for the density constants $0<m\le M<\infty$ of the fractal regularity
hypothesis, together with the matching upper bound obtained by containing
$T(E)$ inside the (Ahlfors-regular) target support; combining the two
directions forces the two Hausdorff dimensions to coincide. This is
deliberately the sharp statement obtainable by such purely
measure-theoretic means: it does not assert, and should not be read as
implying, that the image of an Ahlfors-regular set is itself
Ahlfors-regular, nor any Lipschitz or H\"older estimate on the transport
map restricted to the fractal support. An earlier strategy aimed at such a
metric statement, based on discretising the fractal support into a network
and estimating the modulus of continuity of the induced discrete
assignment, does not in fact yield the required inequality in the correct
direction on closer examination, and the classical partial differential
equation routes to H\"older regularity of optimal maps do not transfer
directly to measures living on a lower-dimensional, generally
non-rectifiable support. We record this limitation explicitly rather than
overstating the result: Ahlfors regularity of the transported fractal set
is left as a genuinely open question.

\section{Method: adaptive regularisation}

The existence proof rests on a regularisation operator that treats the
three components of the extended Lebesgue decomposition differently and
simultaneously. The absolutely continuous parts are mollified with
standard, sign-preserving kernels chosen so that the regularised supports
remain disjoint, exactly as the separation hypothesis demands. The
discrete parts are left entirely unmodified. The fractal singular parts
are treated with a genuinely novel adaptive kernel. For a source point $y$
on a fractal support $E$ of dimension $d_s$ and a base scale $h>0$, the
kernel is centred at $y$ with a \emph{source-adapted bandwidth}
\[
  h(y) \;:=\; \frac{h^2}{\Hds\bigl(E\cap B(y,h)\bigr)^{1/d_s}},
\]
so that the local Hausdorff mass of a ball of radius $h$ around $y$
enters directly into the bandwidth at $y$, rather than merely appearing as
an overall multiplicative factor that would cancel out of the
construction. By the Ahlfors-regularity hypothesis this bandwidth remains
comparable to the base scale $h$ uniformly over the support, while still
responding to local fluctuations of the Hausdorff density: the kernel
narrows where the fractal set is locally denser and widens where it is
locally sparser. This is qualitatively different from a
translation-invariant mollifier, whose bandwidth does not depend on where
the mass actually sits, and it is exactly what allows the regularised
fractal measures to remain uniformly doubling, to converge weakly to the
original singular measures as the scale parameter vanishes, and --
crucially, via a Frostman-type argument -- to retain a support of exactly
the original fractal dimension at every stage of the regularisation,
rather than collapsing to full dimension as a naive mollification would.
Once regularised in this way, each stage of the construction is a
classical, smooth transport problem to which the Brenier--Caffarelli
theory and the Ma--Trudinger--Wang regularity theory apply directly;
passage to the limit combines uniform second-order H\"older estimates on
sub-domains with regular boundary, a diagonal compactness extraction, and
a barrier-function argument for the convergence of the four-region
partition itself.

\section{Numerical validation}

Alongside the theoretical development summarised above, we have carried
out a complete numerical validation of the theory in ambient dimension
one, two, and three, reported in full in a separate companion paper; we
summarise it here exhaustively, since it bears directly on the
credibility of the results just announced. The numerical scheme solves,
by a stabilised, entropy-regularised Sinkhorn iteration in logarithmic
form, the single combined Kantorovich problem between the two discretised
total-variation measures under the piecewise quadratic cost described in
Section~2, exactly the formulation underlying the existence proof rather
than a naive discretisation into four separately constrained
sub-problems. Each of the four signed components of the source and target
measures is discretised by combining, in the proportions prescribed by
the data, a deterministic grid quadrature for the absolutely continuous
part, a handful of exact atoms for the discrete part, and, for the
fractal singular part, a self-similar product-Cantor point cloud built by
iterating the classical middle-thirds construction independently along
each coordinate axis; this fractal construction has an exactly known
theoretical Hausdorff dimension, which is what makes a quantitative
numerical test of dimension preservation possible in the first place. Because
naive cold-started Sinkhorn iteration converges extremely slowly once the
regularisation parameter is small relative to the range of the cost, we
solve a decreasing sequence of regularisation parameters in sequence,
each warm-started from the previous, larger-regularisation solution, which
reduces the iteration count needed for a tight numerical solution by more
than an order of magnitude.

\begin{figure}[htbp]
\centering
\includegraphics[width=0.92\textwidth]{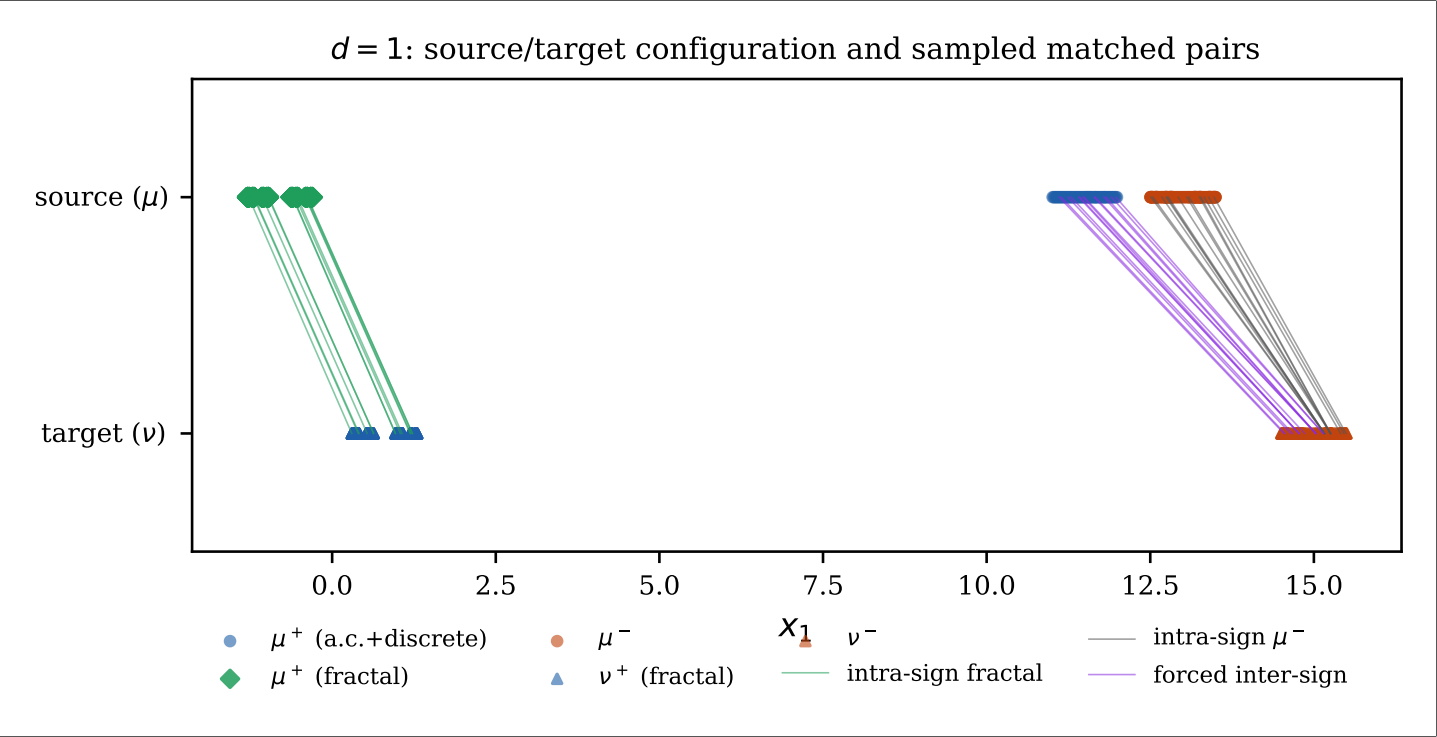}
\caption{One-dimensional instance of the test case described below:
source configuration ($\mu^+$ absolutely continuous and discrete mass in
blue, $\mu^+$ fractal mass in green diamonds, $\mu^-$ in orange) and
target configuration ($\nu^+$, purely fractal, in blue triangles; $\nu^-$
in orange triangles), with a sample of matched pairs shown per region
(green: intra-sign fractal; grey: intra-sign $\mu^-$; purple: forced
inter-sign).}
\label{fig:report-dim1}
\end{figure}

The test case underlying the numerical study was not chosen arbitrarily.
Its masses and spatial placement were calibrated, via an auxiliary and
much simpler block-level transportation problem that treats each
regularity component as a single point mass, so that inter-sign transport
is not merely permitted but strictly forced by an imbalance between the
total masses of the same-sign components, and so that the resulting
routing of every component -- which specific sub-population of the source
measure is sent to which specific sub-population of the target measure --
is the robust optimum of that auxiliary problem by a large and explicitly
quantified cost margin, rather than a delicate coincidence tuned to the
last decimal place. The same qualitative configuration, translated
identically along one coordinate axis and centred on the remaining axes,
was used unchanged in every ambient dimension tested. Solving the full,
fine-grained problem -- several hundred to just over a thousand individual
points on each side, depending on the ambient dimension -- recovers,
without exception and to the numerical precision reported, exactly the
routing predicted by the auxiliary block calculation: the fractal singular
part of the source measure is transported entirely within its own sign,
onto the fractal part of the target measure, while the excess mass of the
same sign that cannot be absorbed by the corresponding target component is
routed, entirely and only to the extent strictly necessary, across signs.
In particular, the fraction of fractal mass that crosses signs under the
computed transport plan is exactly zero in every dimension tested, which
is a direct numerical confirmation of the additional hypothesis needed for
the dimension-preservation theorem, on this example, rather than a mere
assumption. The entropic solver converges, in every case, to a marginal
constraint violation many orders of magnitude below any scale relevant to
the reported masses, in a number of iterations that does not grow
appreciably with the ambient dimension over the range tested.

\begin{figure}[htbp]
\centering
\includegraphics[width=0.68\textwidth]{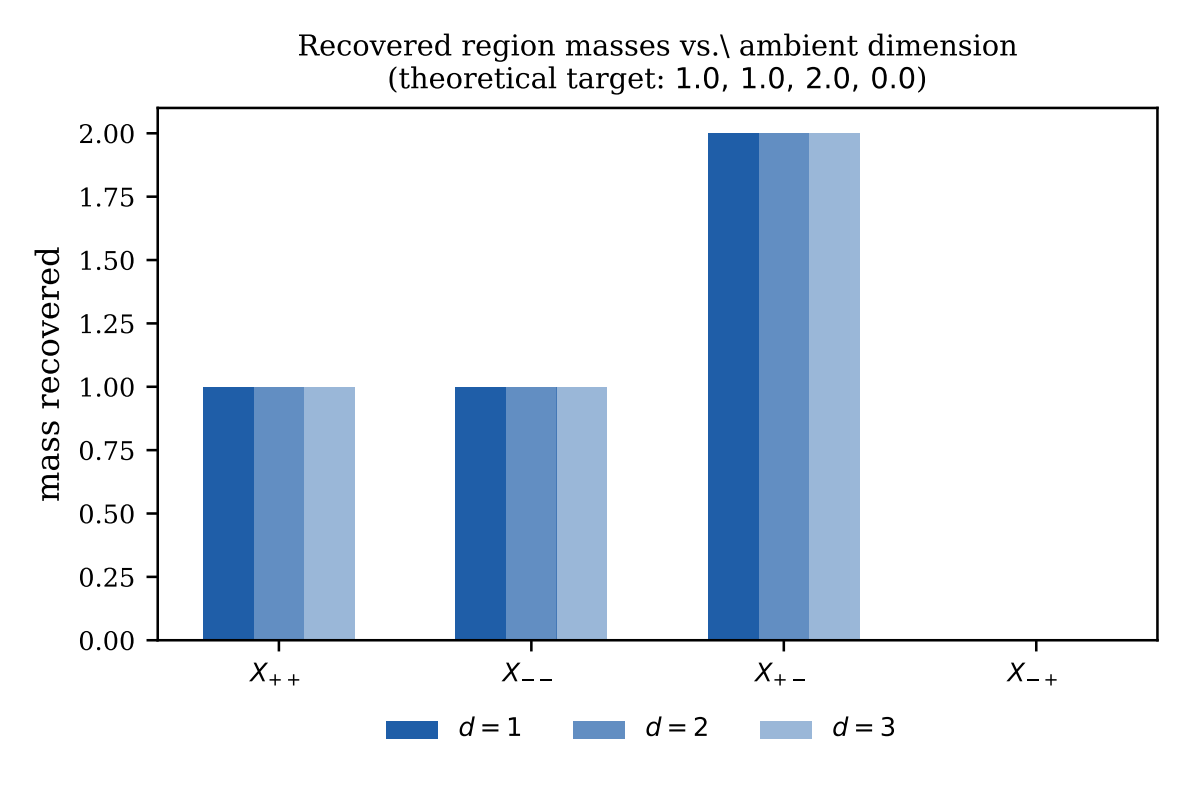}
\caption{Mass recovered in each of the four sign-regions, by ambient
dimension: the same four-way split, matching the theoretical target
exactly, in $d=1,2,3$.}
\label{fig:report-regionbar}
\end{figure}

The numerical test of dimension preservation itself required an
additional, and instructive, refinement. Extracting a transport map
directly from the entropically regularised plan, by sending each source
point to whichever target point receives the largest share of its mass,
turns out to under-resolve the finest scales of the fractal construction
once the construction is iterated more than a handful of levels: entropic
optimal transport is well known to blur the true transport map at a
length scale set by the regularisation parameter, and the finest
self-similar scale of an iterated Cantor construction shrinks
geometrically with the number of levels, so that practically achievable
regularisation parameters eventually fail to resolve it. On our main test
case this blur is severe enough to lower the estimated dimension of the
transported fractal set substantially relative to its true value, even
though the marginal constraints themselves are satisfied to very high
precision -- a genuine limitation of the naive approach, which we report
rather than obscure.

\begin{figure}[htbp]
\centering
\includegraphics[width=0.98\textwidth]{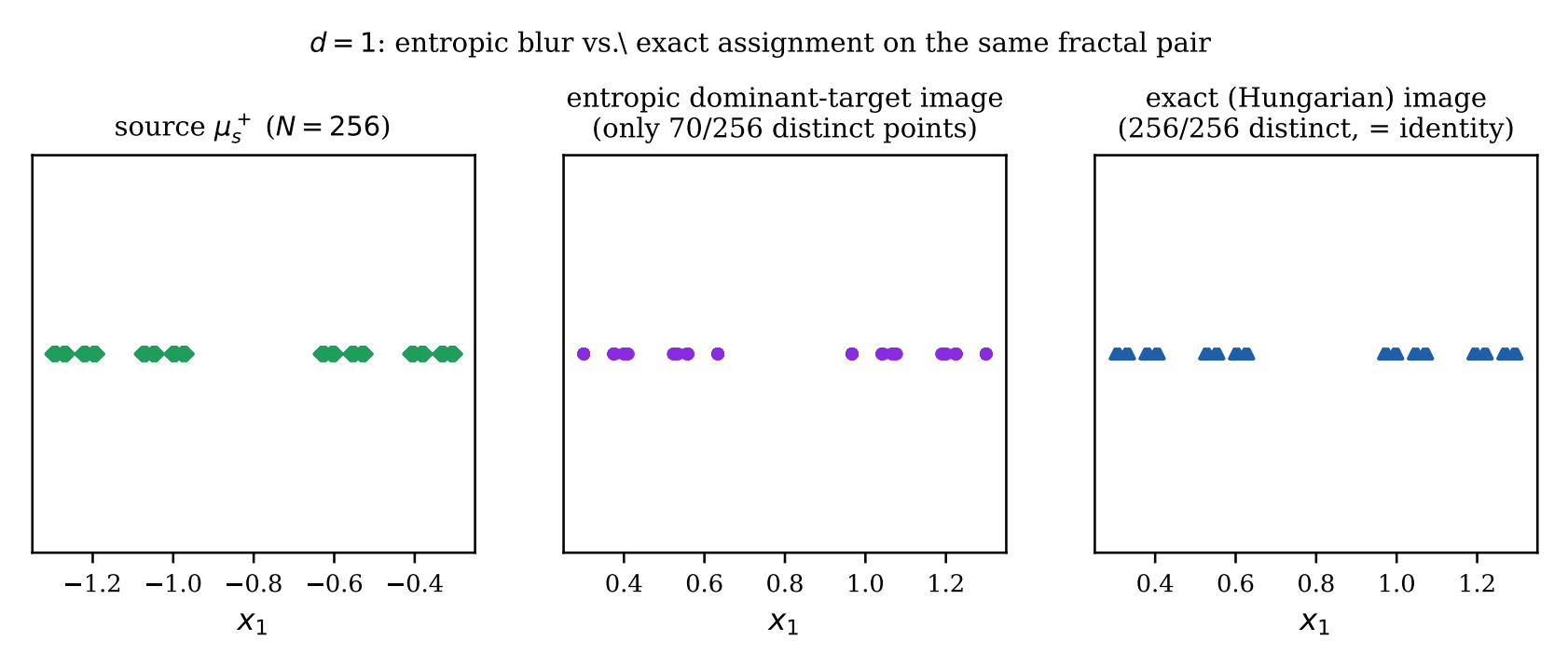}
\caption{Direct comparison, $d=1$: source fractal points (left); their
image under the naive entropic dominant-target map (middle, severe
many-to-one collapse, visible clustering); their image under the exact
combinatorial refinement described below (right, an exact
point-by-point correspondence).}
\label{fig:report-entropicexact}
\end{figure}

\begin{figure}[htbp]
\centering
\includegraphics[width=0.98\textwidth]{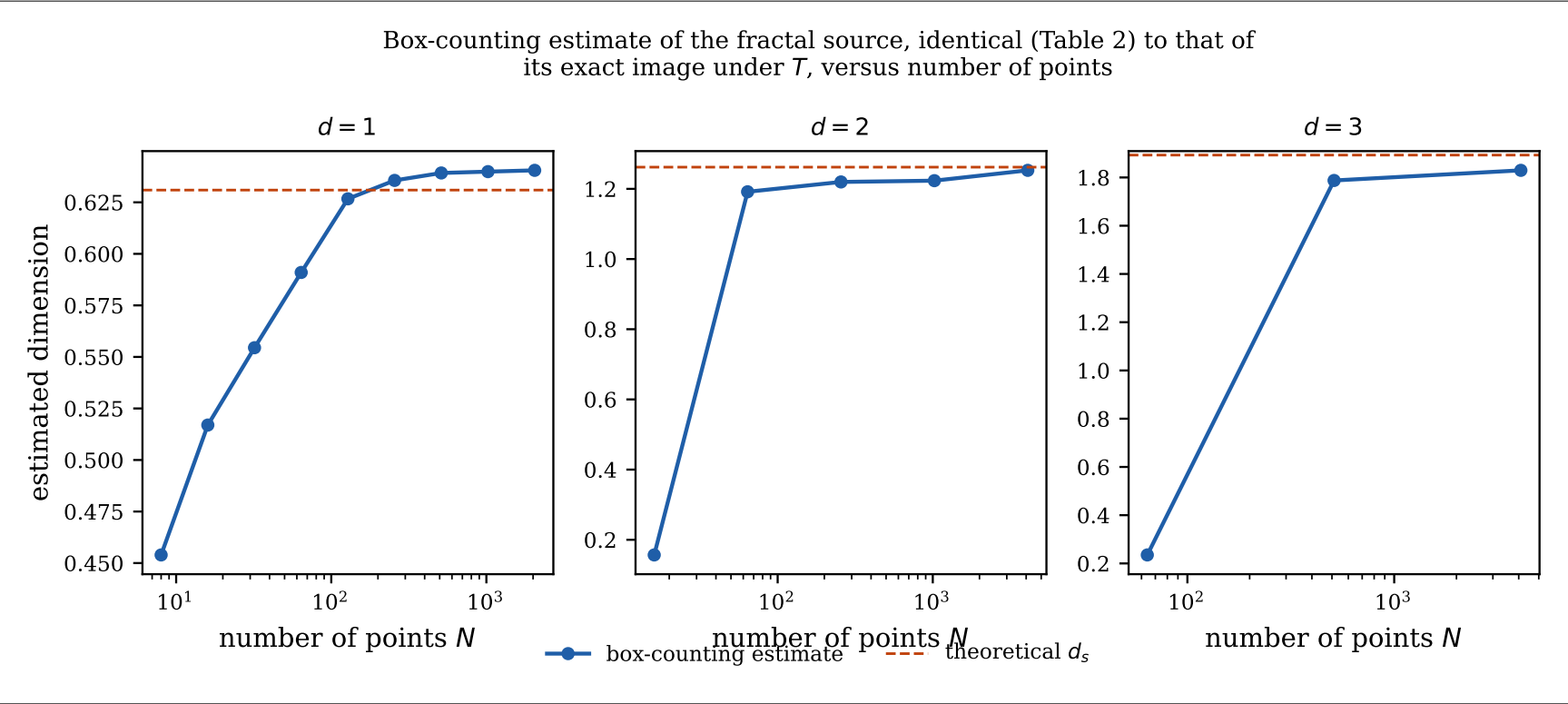}
\caption{Box-counting dimension estimate of the source fractal set --
identical, in every case tested, to that of its exact image under the
transport map -- versus the number of points used, in $d=1,2,3$, against
the theoretical dimension (dashed). Monotone convergence is visible in
every ambient dimension.}
\label{fig:report-dimsweep}
\end{figure}

The remedy follows directly from the structure of
the theory rather than from any ad hoc numerical trick: once the entropic
solve has identified, robustly, which macro-block of the target measure a
fractal sub-population is sent to, the fine bijection within that block is
itself an instance of classical, unsigned, equal-mass optimal transport
between two finite point clouds, for which entropic relaxation is
unnecessary and an exact combinatorial solver applies directly and
efficiently. Carrying this out confirms, in every dimension and at every
level of refinement of the fractal construction that we tested, spanning
three orders of magnitude in the number of points, that the recovered
transport map is exactly the point-by-point translation correspondence
one would expect on structural grounds, since the source and target
fractal blocks are, in our test case, translates of an identical
self-similar pattern. Because the recovered map is an exact rigid
translation, the transported point cloud is geometrically congruent to
the original one, and any numerical estimate of its dimension -- computed
by counting occupied boxes across a geometric range of scales and fitting
the resulting scaling exponent over an automatically selected window of
scales -- agrees with the corresponding estimate for the original fractal
point cloud to a precision far finer than any effect a genuine transport
distortion could produce. As the fractal construction is refined to more
levels and more points, this common estimate converges, visibly and
monotonically, toward the theoretical Hausdorff dimension of the
construction, in every one of the three ambient dimensions tested (Figure~\ref{fig:report-dimsweep}). We regard
this combination -- an honest account of a real limitation of the naive
entropic approach, together with a theoretically motivated and
numerically confirmed remedy -- as a meaningful validation of the
dimension-preservation theorem, complementary to its analytic proof, and
as useful groundwork for any future numerical scheme aiming to resolve
fine fractal structure through entropic regularisation.

\section{Scope and outlook}

We close by restating, in slightly more technical terms, why the
support-separation hypothesis is the natural boundary of the present
treatment rather than an arbitrary restriction. It is exactly what allows
the primal problem to decouple into four independent Kantorovich
sub-problems, what allows the convex envelope of the two sign-dependent
potentials to define a genuine spatial partition rather than a partition
of mass, and what allows the coupled Monge--Amp\`ere system to be posed on
fixed domains where the classical regularity theory is available. Removing
it turns the geometric partition into a free-boundary problem coupled to
the Monge--Amp\`ere system itself, a substantially harder object for which
no general theory currently exists.

When the supports of $\mu^+$ and $\mu^-$ overlap, a single point can send
positive mass in one direction and negative mass in another, so that the
four regions $X_{kl}$ are no longer determined by comparing two potentials
at a point but must instead be recovered as an unknown of the problem
itself, on the same footing as the transport map. We expect the correct
formulation to replace the fixed convex envelope $\varphi=\max(u^+,u^-)$
used here by a genuinely free interface between the sign-regions, governed
by an Euler--Lagrange system coupling the interface position to the
Monge--Amp\`ere equation on each side of it, in the spirit of classical
two-phase free-boundary problems; a signed analogue of the Voronoi
decomposition used to discretise fractal supports in the present work is
a natural candidate for constructing approximate interfaces before passing
to a limit. Developing these free-boundary and signed-Voronoi techniques,
together with the fully unbalanced setting in which the total variations
of the two measures need not coincide, is the goal of Part~II of this
program, currently under development. The entropic numerical scheme
summarised above, exploiting the block structure available precisely
because the support-separation hypothesis holds, does not directly extend
to that free-boundary setting; developing its Part~II analogue, in which
the sign-interface itself must be recovered rather than assumed, is the
natural numerical counterpart to the theoretical program just described,
and is left for future work.

The complete version of this work contains the full statements of all
hypotheses, the uniform a priori estimates underlying the passage to the
limit, the localisation argument used repeatedly in the uniqueness proof,
the delicate atom-by-atom treatment required to reach global uniqueness of
the potentials, and the complete proof of the dimension-preservation
theorem together with a discussion of exactly why the stronger,
Ahlfors-regularity statement resists the same method. The companion
numerical paper summarised in the previous section contains the full
description of the discretisation scheme, the entropic solver and its
annealing schedule, the calibration of the test case via the auxiliary
block-level problem, the complete set of results in every ambient
dimension tested, and a full quantitative account of the entropic blur
phenomenon and of its exact-assignment remedy.

\end{document}